\documentclass[12pt]{amsart}
\usepackage{amscd}
\advance\hoffset by -0.5 truecm
\usepackage{amssymb}
\newtheorem{Theorem}{Theorem}[section]
\newtheorem{Lemma}[Theorem]{Lemma}
\newtheorem{Corollary}[Theorem]{Corollary}

\newtheorem{Proposition}[Theorem]{Proposition}

\newtheorem{Remark}[Theorem]{Remark}

\def \dim{{\mbox {dim}}\,}

\def\V{\mbox{Var}}

\def\Z{{\mathbb Z}}
\def\R\re
\def\V{\bf V}

\def \la{\lambda}

\def \re{{\mathbb R}}
\def \Q{{\mathbb Q}}
\def \cp{{\mathbb CP}}
\def \T{{\mathbb T}}
\def \C{{\mathbb C}}
\def \M{{\widetilde{M}}}

\def \0{\lambda_{0}}
\def \la{\lambda}

\def\h{{\rm h}_{\rm top}(g)}

\begin{document}
\title{Zero entropy and bounded topology}

\author[G. P. Paternain]{Gabriel P. Paternain}\thanks{G.P. Paternain
was partially supported by CIMAT, Guanajuato, M\'exico.}
\address{Department of Pure Mathematics and Mathematical Statistics,
University of Cambridge,
Cambridge CB3 0WB, England}
\email{g.p.paternain@dpmms.cam.ac.uk}

\author[J. Petean]{Jimmy Petean}
 \address{CIMAT  \\
          A.P. 402, 36000 \\
          Guanajuato. Gto. \\
          M\'exico.}
\email{jimmy@cimat.mx}

\thanks{J. Petean is supported by grant 37558-E of CONACYT}



\begin{abstract}
We study the existence of Riemannian metrics with zero topological
entropy on a closed manifold $M$ with infinite fundamental group. We show 
that such a metric does not exist if there is a 
finite simply connected CW
complex  which maps to $M$ in such a way that the rank of the map
induced in the pointed loop space homology grows exponentially. This
result allows us to prove in dimensions four and  five,  
that if $M$ admits a metric with zero
entropy then its universal covering has the rational homotopy type
of a finite elliptic CW complex. We conjecture that this is the case 
in every dimension.

\end{abstract}

\maketitle

\section{Introduction}

Let $M^{n}$ be a closed connected smooth manifold.
Given a Riemannian metric $g$, let $\phi_{t}$ be the geodesic flow
of $g$.

One of the most fundamental dynamical invariants that one can associate to
$\phi_{t}$ is the {\it topological entropy,} which we denote by $\h$. 
It roughly measures the orbit structure complexity of the flow.
Positive entropy means in general, that the geodesic flow presents
somewhere in the phase space (the unit sphere bundle of the
manifold) a complicated dynamical behaviour. There are various
equivalent ways of defining entropy, but for the geodesic flow, Ma\~n\'e's formula
\cite{Man} provides a clear understanding of this invariant in terms of geodesic arcs.
Given points $p$ and $q$ in $M$ and $T>0$, define $n_{T}(p,q)$ to be the number of
geodesic arcs joining $p$ and $q$ with length $\leq T$. We have
 $$\h= \lim_{T\rightarrow
\infty}\frac{1}{T}\log \int_{M\times M}n_{T}(p,q)\;dp\,dq. $$

The main goal of this paper is to address the following natural question:
which manifolds admit metrics with zero topological entropy?

A classical result of E.I. Dinaburg \cite{D} asserts that if $M$ admits such a metric,
then $\pi_1(M)$ must have subexponential growth. It is still unknown if there are
finitely presented groups which are of subexponential growth, but not of polynomial growth.
If such groups do not exist, then zero entropy implies that $\pi_1(M)$ is virtually
nilpotent, thanks to a celebrated theorem of M. Gromov \cite{G}.
For closed geometrizable 3-manifolds, this obstruction on the fundamental group
is enough to determine those which admit a metric with zero topological entropy, cf. \cite{AP}.

In the late 1980's new topological obstructions were found, this time for simply connected
manifolds. Y. Yomdin \cite{Y} proved a fundamental theorem for general $C^{\infty}$
dynamical systems relating the topological entropy with the volume growth of submanifolds
which paved the way to Ma\~n\'e's formula.
When combined with the Morse theory of the loop space and a beautiful discovery of
Gromov \cite{G1} concerning cycles with bounded length in the pointed loop space $\Omega M$,
it gave strong restrictions to zero entropy. Namely, if $M$ is simply connected and admits a $C^{\infty}$ metric $g$ with $\h=0$, then the sum of the Betti numbers
$\sum_{i=1}^{n}\mbox{\rm dim}\,H_{i}(\Omega M,k_{p})$ grows subexponentially
with $n$ for any field of coefficients $k_p$, $p$ prime or zero.
When $k_p=\Q$, this implies that $M$ is rationally elliptic, i.e. $\pi_*(M)\otimes\Q$
is finite dimensional (cf. \cite{FHT}). We refer to \cite{P} for an account of these developments.

However, these results only hold for simply connected
manifolds (or finite $\pi_1(M)$) because Gromov's theorem does require to control
the length of paths running on the 1-skeleton of a triangulation of $M$ and in the simply
connected case we can always collapse the 1-skeleton to a point by a map homotopic to the
identity which gives the desired control.

What topological restrictions to zero entropy do we have when $\pi_1(M)$ is infinite and 
of subexponential growth? We begun looking at this problem in \cite{PP} motivated
by the minimal entropy problem for compact complex surfaces.
Here we show:

\medskip

\noindent {\bf Technical Lemma.} {\it Let $M$ be a closed manifold. 
Let $f:K\to M$ be a continuous map, where $K$ is a finite simply
connected  CW complex
and let $\Omega(f)$ be the induced map between pointed loop spaces. Let
$H_{*}(\Omega(f),k_p)$ be the map induced in homology with some field 
of coefficients $k_p$ and let $R_i$ be the rank of this map in dimension
$i$. Set
$$R:= \limsup_{i\to \infty}\frac{1}{i}\log\left(\sum _{j\leq i} R_{j}\right).$$ 
If $R>0$, then given any smooth Riemannian metric $g$ on $M$ we have:
\[\h > \frac{\la(g)}{2}.\]

}

\medskip

In the inequality, $\la(g)$ is the {\it volume entropy} 
of the Riemannian manifold which
is defined as the exponential growth rate of the volume of balls 
in the universal covering of $M$.
Recall that Manning's inequality
\cite{Ma} asserts that for any metric $g$,
${\rm h}_{\rm top}(g)\geq \la(g)$ and it is well known that 
$\la(g)>0$ if and only if $\pi_1(M)$ has exponential
growth. We are interested in the inequality as an obstruction to
the existence of metrics with vanishing topological entropy, particularly
in the case when $\pi_1(M)$ has subexponential growth. It would be quite
interesting to be able to replace $\la(g)/2$ by $\la(g)$ in the 
Technical Lemma.

It seems useful to note the following point: if $\M$ is the
universal covering of $M$, then the projection induces an isomorphism
between the homology of the loop space of $\M$ and the
homology of the connected component of the loop space of $M$ given by the contractible loops. 
Therefore one can consider a CW complex $K$ which maps to $\M$ and then compose with the 
projection to $M$ to be in the conditions of the Technical Theorem. We will
use this remark in all of our examples.

Recall that a connected CW complex $X$ is said to be {\it nilpotent} if $\pi_1(X)$ is a nilpotent group and operates nilpotently on $\pi_i(X)$ for 
every $i\geq 2$.
As an immediate corollary of the lemma we have:

\medskip

\noindent {\bf Corollary.} {\it Let $M$ be a closed nilpotent manifold.
If $M$ admits a smooth metric with zero topological entropy, then
$\pi_{*}(\Omega M)\otimes \Q$ is finite dimensional.}

\medskip

Indeed, if $M$ is nilpotent, all the homology groups of $\M$ are finitely 
generated (cf. \cite[Theorem 2.16]{HMR}) and thus there is a finite simply
 connected CW complex $K$ and a homotopy equivalence $f:K\to \M$.
The complex $K$ must be rationally elliptic by the Technical Lemma.

We proved the lemma in \cite[Theorem C]{PP} when $K$ is a smooth
compact  manifold
with boundary which is embedded in $\M$ and for which the
corresponding map in the loop space homology is an injection. 
The disadvantage of this earlier version
is that in order to use it we need to have some apriori knowledge of $\M$
so that we can find our embedded $K$, while with the current version
$K$ and $f$ will arise by simple topological considerations as in the corollary above.

Nevertheless the old version was good enough to prove results like the 
following \cite[Theorem D]{PP}: if $M$ admits a metric with zero entropy and it
can be decomposed as $X_{1}\# X_{2}$, where the order of
the fundamental group of $X_1$ is at least $3$, then $X_2$ is a homotopy sphere.

We now pose the main topological question that the Technical Lemma suggests:

\medskip

\noindent {\bf Question.} {\it Let $M$ be a closed manifold whose fundamental group
has subexponential growth. 
If $\mbox{\rm dim} \,H_{*}(\M,\Q)=\infty$, does there exist
a finite 1-connected rationally hyperbolic complex
$K$ and a map $f:K\to M$ for which the rank of $H_{*}(\Omega(f),\Q)$ 
grows exponentially? }

\medskip

Of course one can formulate similar questions for other fields of coefficients, but we 
believe it should be easier to deal first with the case of characteristic zero, due to the
technology at our disposal provided by Rational Homotopy Theory.

If the Question has a positive answer, then the Technical Lemma implies that if $M$
admits a metric with zero entropy then $\M$ has the rational homotopy type
of a finite elliptic 1-complex, so we see that zero entropy implies bounded
topology in a very strong sense.

In the present paper we will prove by simple topological arguments that
the Question has a positive answer in dimensions 4 and 5.
This in turn will give us an essentially complete picture of which 4-manifolds
have metrics with zero entropy and will allow us to close some gaps left open
in \cite{PP}.

Let us describe these results in more detail. 
From now on if in ordinary homology coefficients are not indicated they
are meant to be $\Z$. In the next theorem, $\sigma$ and $\chi$ stand
for signature and Euler characteristic respectively.

\medskip
\noindent {\bf Theorem A.} {\it
Let $M$ be a closed 4-manifold with infinite fundamental group $\pi$.
If $M$ admits a metric with zero topological entropy, then $\sigma (M)=\chi(M)=0$,
$\M$ has the rational homotopy type of a finite simply connected 
elliptic CW complex and $H_{2}(\M)\cong H^2(\pi,\Z[\pi])$.
Moreover, if we assume further that $\pi$ has polynomial growth then, $M$ is
finitely covered by one the following:
\begin{enumerate}
\item $S^3\times S^1$;
\item a manifold s-cobordant to $S^2\times \T^2$;
\item a manifold homeomorphic to a nilmanifold.
\end{enumerate}}

\medskip

Using the results on 4-manifolds, we can now complete the classification
of compact complex surfaces which admit a metric with
zero entropy. We begun this classification in \cite{PP}, but our results
excluded two cases: surfaces of general type and
surfaces of class VII with positive second Betti number.
It is unknown if there are surfaces of general type homeomorphic
to $S^{2}\times S^{2}$ or $\cp^{2}\#\overline{\cp}^{2}$, although
it is known that there is no surface of general type {\it diffeomorphic}
to $S^{2}\times S^{2}$ or $\cp^{2}\#\overline{\cp}^{2}$.
We call such a potential example, {\it an exotic surface of general type}.
In the next theorem we view compact complex surfaces as smooth 4-manifolds
and we ignore their complex structures.

\medskip

\noindent {\bf Theorem B.} {\it Let $S$ be a compact complex surface which is
not an exotic surface of general type.
Then $S$ admits a metric with zero topological entropy if and only if
$S$ is diffeomorphic to one of the following:
$\C P^2$, a ruled surface of genus $0$ or $1$, a complex torus, a hyperelliptic surface,
a Hopf surface, a Kodaira surface, or a Kodaira surface modulo a finite group.}

\medskip

Finally in dimension 5 we prove:

\medskip

\noindent {\bf Theorem C.} {\it Let $M$ be a closed 5-manifold with infinite fundamental group. 
If $M$ admits a metric
with zero entropy, then $\M$ has the rational homotopy type of a finite 1-connected
elliptic complex. Moreover, $H_{3}(\M)\cong H^2(\pi,\Z[\pi])$.}

\medskip

In fact, for most groups $\pi$ with subexponential growth the second end
group $H^2(\pi,\Z[\pi])$ is either $0$ or $\Z$. Our methods also yield information
at the torsion level. For example we will show that if $H^2(\pi,\Z[\pi])$ is either
$0$ or $\Z$, then $H_{2}(\M)$ has no finite subgroup as a direct summand.

\medskip

{\it Acknowledgements:} We thank Burt Totaro for several useful
comments on the first draft of the manuscript.

\section{Proof of the Technical Lemma}

Let $(M,g)$ be a Riemannian manifold and let $K$ be a finite simply
connected CW complex. Given a continuous map $f:K\rightarrow M$ we let
$\Omega (f) :\Omega (K) \rightarrow \Omega (M)$ 
be the obvious map induced between the corresponding pointed
loop spaces. The following lemma is essentially due to Gromov \cite{G1,bates}:

\begin{Lemma} There exists a constant $C=C(K,f,M,g)$ such
that given any homology class $\psi \in H_i (\Omega (K))$, the
class $\Omega (f)_* (\psi )$ can be represented by a cycle of
Lipschitz curves in $M$ with length bounded by $Ci$.
\end{Lemma}

\begin{proof} Since $K$ is homotopy equivalent to a finite simply
connected simplicial complex \cite[Theorem 2.C5]{hatcher}, we can 
assume that $K$ is actually a simplicial complex. We can consider $K$
as a subcomplex of a simplex $\Delta^N$ and restrict the standard
metric on $\Delta^N$ to give a metric on $K$. Of course, this 
metric restricts to the standard Euclidean metric on each simplex
of $K$. It is easy to see that one can approximate $f$ by a homotopic
map which is Lipschitz with respect to this metric. Therefore we will
also assume that the map $f$ is Lipschitz.

\vspace{.5cm}

Let $h:L\rightarrow \Omega (K)$ be a map from a finite simplicial 
complex of dimension $i$  with an $i$-th homology class mapping to
$\psi$. The map $h$ corresponds to a map $H:L\times [0,1]
\rightarrow K$. Given a positive integer $k$ we will consider the
simplicial structure on $[0,1]$ obtained by subdividing the
interval into $k$ subintervals of equal length. The simplicial 
structures on $L$ and $[0,1]$ give a natural cellular decomposition
on $L\times [0,1]$.

For this cellular decomposition one can obtain a simplicial
approximation similar to the simplicial case as follows:

\vspace{.5cm}

Let $S$ be any simplicial complex. We call a map $R:L\times [0,1]
\rightarrow S$ {\it simplicial} if it is a simplicial map when
restricted to each $L\times \{j/k\}$ and for any $x\in L$, 
$R(x,j/k)$ and $R(x,(j+1)/k)$ belong to a simplex in $S$ and the
restriction of $S$ to the vertical segment $\{ x\} \times
[j/k,(j+1)/k]$ is linear.

Given a continuous map $r:L\times [0,1]\rightarrow S$, we say that
$R$ is {\it a simplicial approximation} of $r$ if it is a simplicial map
such that $R(q)\in \mbox{\rm Carrier}(r(q))$ (the smallest simplex
containing $r(q)$) for any $q\in L\times [0,1]$. It
is easy to see as in the simplicial case that if $R$ is a simplicial
approximation of $r$, then $R$ and $r$ are homotopic.

Recall now that the open star of a vertex in a simplicial complex
is the union of the interior of all the simplices containing the
vertex, and
for a vertex $(v,j/k) \in L \times [0,1]$ define its open star as
$\mbox{\rm Star}(v,j/k)=\mbox{\rm Star}(v) \times ((j-1)/k,(j+1)/k)$. It is clear that
after enough subdivisions of the simplicial structure on $L$ and 
taking $k$ big enough, we can assume that the diameter of the open
star of any vertex is as small as we want. Therefore we can assume 
that for any vertex $(v,j/k)\in L \times [0,1]$  there exists a vertex  
$w \in S$ such that the open star of $(v,j/k)$
is contained in $r^{-1} (\mbox{\rm Star}(w))$. We define $R(v,j/k)=w$. As in the
simplicial case, we can extend $R$ to each $L\times \{ j/k\}$ as a
simplicial approximation of the restriction of $r$ to 
$L\times \{ j/k\}$.

Given $x\in L$, let $\sigma =\mbox{\rm Carrier}(x)$ and let $v_1 ,...,v_l$ be the
vertices of $\sigma$. Let $w_i =R(v_i ,j/k)$ and $w_{i+l} =
R(v_i ,(j+1)/k)$. Note that $\mbox{\rm Interior}(\sigma ) \times (j/k,(j+1)/k)$
is contained in the open stars of each $(v_i ,m/k) (m=j$ or $j+1)$. 
By construction we have that $r(\mbox{\rm Interior}(\sigma )\times
(j/k,(j+1)/k))$ is contained in $\mbox{\rm Star}(w_i )$ for each $i=1,...,2l$.
The intersections of these sets is nonempty and this implies that
those vertices form a simplex in $S$. We have done this to show that
there exists a simplex of $S$ which contains both $R(x,j/k)$ and
$R(x,(j+1)/k)$ and so we can extend $R$ as a simplicial map on
$L\times [0,1]$.

Let us finally check that $R$ is actually a simplicial approximation
of $r$. We have to show that for any $(x,t)\in L\times [0,1]$,
$R(x,t)\in \mbox{\rm Carrier}(r(x,t))$. We already know this if $t$ is a vertex
of $[0,1]$. So we can assume that $t\in (j/k,(j+1)/k)$ for some $j$.
As in the previous paragraph, let $v_1 ,...,v_l$ be the vertices of
$\mbox{\rm Carrier}(x)$ and let $w_1 ,...,w_{2l}$ be the corresponding vertices
in $S$. We have that $(x,t)\in \mbox{\rm Star}(v_i ,j/k) \cap \mbox{\rm Star} 
(v_i ,(j+1)/k)$ for each $i=1,...,l$. Then 
$r(x,t)\in r(\mbox{\rm Star}(v_i ,j/k) 
\cap \mbox{\rm Star}(v_i ,(j+1)/k)) \subset \mbox{\rm Star}(w_i )\cap \mbox{\rm Star}(w_{i+l})$. 
This implies that $w_i ,w_{i+l} \in \mbox{\rm Carrier}(r(x,t))$. Since $R(x,t)$
is a linear combination of $w_1 ,...,w_{2l}$ and all these vertices
belong to the simplex $\mbox{\rm Carrier}(r(x,t))$, we get that $R(x,t)\in
\mbox{\rm Carrier}(r(x,t))$ and therefore $R$ is a simplicial approximation of
$r$.

\vspace{.8cm}

Therefore we can take a simplicial approximation $R$ of the map
$H:L\times [0,1] \rightarrow K$ and consider the corresponding map
$r:L\rightarrow \Omega (K)$. Let us consider the space 
${\Omega (K)}_{pl,k} \subset \Omega (K)$ given by those paths
which are linear on each segment of the form $[j/k,(j+1)/k]$. Note
that by construction there exists a $k$ such that the image of
$r$ is contained in ${\Omega (K)}_{pl,k}$. Each element in
${\Omega (K)}_{pl,k}$ determines a point in $K^{k-1}$ and in this way we 
identify ${\Omega (K)}_{pl,k}$ with a subset of $K^{k-1}$ (recall that
the initial and final points are fixed). The simplicial structure on
$K$ induces a cellular decomposition in $K^{k-1}$ and 
${\Omega (K)}_{pl,k}$ is a subcomplex: it is the union of all the
cells $\sigma_1 \times ...\times {\sigma}_{k-1}$ such that 
$\sigma_i $ and ${\sigma}_{i+1}$ are contained in a simplex of
$K$ and the same for the initial point and $\sigma_1$ and the 
end point and ${\sigma}_{k-1}$. 

After another homotopy we can assume that the image of $r$ is
contained in the $i$-th skeleton of ${\Omega (K)}_{pl,k}$ with
respect to the cell decomposition described above.

Since $L$ is simply connected there exists a simplicial map
$\alpha :L\rightarrow L$ homotopic to the identity and which maps
the whole 1-skeleton of $L$ to a point.

Now let $c$ be a path in the $i$-th skeleton of  ${\Omega
  (K)}_{pl,k}$. This means that $c$ belongs to a cell of the form 
$\sigma_1 \times ...\times {\sigma}_{k-1}$ with $\Sigma_j \mbox{\rm dim}(\sigma_j
  ) \leq i$. Thus, the path $c$ is formed by segments joining a
pair of vertices of the triangulation and at most $2i$ segments in
which one of the points is not a vertex. After composing with the map
$\alpha$ the former are sent to a point while the latter are sent
to paths of length bounded by a constant which depends on $\alpha$
but not on $i$. Hence there exists a constant $C'$ such that the
image of the $i$-th skeleton of  ${\Omega (K)}_{pl,k}$ is sent by
composition with $\alpha$ to a set of paths with length bounded
by $C' i$. In this way we see that we can represent $\psi$ by
a cycle formed with paths with length bounded by $C' i$. Composing
with $f$ we see that we can represent $\Omega (f)_* (\psi )$ by
a cycle formed with paths  of length bounded by $Ci$, where $C$ is
a constant dedending only on $\alpha$ and the Lipschitz constant
of $f$.

\end{proof}

We will use the lemma in the form of the following corollary:

\begin{Corollary} Let $(N,g)$ be a connected complete Riemannian
manifold. Let $K$ be a finite simply connected CW complex  and 
$f:K\rightarrow N$ be a continuous map. 
Denote by $H_i (\Omega(f),k_p )$ the induced map between  
the $i$-th homology groups of the corresponding pointed loop spaces 
(for some field of  coefficients $k_p$) and let $R_i$ be the rank of this map.
Then, there exists a positive constant
$C$ depending only on $K,f,g$ such that for any $x\in K$, $T\geq C i$
and any $y\in B(f(x), T/2)$ we have that $b_i (\Omega^T (N,f(x),y),k_p )
\geq R_i$.
\end{Corollary}

\begin{proof} We know from the lemma that for any $i$-th homology class
$\psi$ in $\Omega (K,x,x)$, $\Omega (f)_* (\psi )$ can be represented 
by a cycle in
$\Omega^{C'i} (N, f(x),f(x))$. Consider now a minimizing geodesic 
between $f(x)$ and $y$. Following the paths in the cycle by this
geodesic we obtain a cycle in $\Omega^{C'i + d(f(x),y)}(N,f(x),y)$.
If $C=2C'$, $T\geq Ci$ and $y\in B(f(x),T/2)$ we get that our new
cycle is in $\Omega^T (N,f(x),y)$ and the corollary follows.
\end{proof}

We are now ready to prove the Technical Lemma in the introduction:

\medskip

\noindent {\bf Technical Lemma.} {\it Let $M$ be a closed manifold. 
Let $f:K\to M$ be a continuous map, where $K$ is a finite simply
connected  CW complex
and let $\Omega(f)$ be the induced map between pointed loop spaces. Let
$H_{*}(\Omega(f),k_p)$ be the map induced in homology with some field 
of coefficients $k_p$ and let $R_i$ be the rank of this map in dimension
$i$. Set
$$R:= \limsup_{i\to \infty}\frac{1}{i}\log\left(\sum _{j\leq i} R_{j}\right).$$ 
If $R>0$, then given any smooth Riemannian metric $g$ on $M$ we have:
\[\h > \frac{\la(g)}{2}.\]

}

\medskip

\begin{proof} By the lifting property of covering spaces we
can assume that we have a map $f:K\to \M$.
Let us recall that for any  $x\in M$ (cf. \cite{P})
\[\h\geq
\limsup_{T\to\infty}\frac{1}{T}\log\int_{M}n_{T}(x,y)\,dy.\]

Let $p:\M\to M$ be the covering projection. It is easy to check that
given any $x\in \M$ we have
\[\int_{M}n_{T}(p(x),y)\,dy=\int_{\M}n_{T}(x,y)\,dy=\int_{B(x,T)}n_{T}(x,y)\,dy.\]
Thus for any $x\in \M$ we have
\begin{equation}
\h\geq
\limsup_{T\to\infty}\frac{1}{T}\log\int_{B(x,T)}n_{T}(x,y)\,dy.
\label{equis} \end{equation}

Now assume that $x=f(z)$ for some $z\in K$. Morse theory tells us
that if $y$ and $x$ are not conjugate, then 

$$n_T (x,y)\geq \sum_{j\geq 0} b_j (\Omega^T (\M,x,y),k_p ).$$
But now, using the previous corollary we get that if $y\in B(x,Ci/2)$,

$$ \sum_{j\geq 0} b_j (\Omega^{Ci} (\M,x,y),k_p ) \geq \sum_{j\leq i}
R_j ,$$

\noindent
where $C$ is the constant appearing in the corollary.
Integrating the previous inequalities with respect to $y\in B(x,Ci/2)$
yields:

$$\int_{B(x,Ci/2)} n_{Ci} (x,y)\,dy \geq \left(\sum_{j\leq i} R_j
\right) \mbox{\rm Vol}(B(x,Ci/2)).$$
Therefore we obtain:

\begin{align*}
\h &\geq 
\limsup_{i\to\infty}\frac{1}{Ci}\log  \int_{B(x,Ci)}n_{Ci}(x,y)\,dy\\
&\geq \limsup_{i\to\infty}\frac{1}{Ci}\log \left(\sum_{j\leq i} R_j \right)
\mbox{\rm Vol}(B(x,Ci/2))  .\\
\end{align*}
And thus, if $R$ is the exponential growth rate of $\sum_{j\leq i} R_j$,

$$\h \geq \frac{R}{C} \, +\frac{\lambda(g)}{2}.$$

\end{proof}

\section{Topological preliminaries}

In the following sections we will try to apply the Technical
Lemma to find obstructions to zero entropy. This is of course  a
purely topological problem and in this section we will summarize
some general techniques and concepts we will use.

\subsection{Domination} Recall that a topological space $Z$ is 
{\it dominated} by
 a topological space $Y$ if
there exist continuous maps 
$r:Y\to Z$ and $\iota:Z\to Y$ such that $r\,\iota$ is homotopic
to the identity of $Z$. 
In several occasions we will make use of the following lemma, sometimes without explicit mention
to it.

\begin{Lemma} Let $X$ and $Y$ be simply connected spaces which are
rational homotopy equivalent. Suppose $Z$ is dominated by $Y$.
Then there is a map $g: Z\to X$ such that
$$H_{*}(\Omega(g),\Q):H_{*}(\Omega(Z),\Q)\to H_{*}(\Omega(X),\Q)$$
is an injection.
\label{domrat}
\end{Lemma}

\begin{proof} Let $f:Y\to X$ be a rational homotopy equivalence.
Since $Y$ dominates $Z$, there is a map $\iota:Z\to Y$ such that
$H_{*}(\Omega(\iota),\Q):H_{*}(\Omega(Z),\Q)\to H_{*}(\Omega(Y),\Q)$
is an injection.
By the Whitehead-Serre theorem, cf. \cite[Theorem 8.6]{FHT},
$H_{*}(\Omega(f),\Q):H_{*}(\Omega(Y),\Q)\to H_{*}(\Omega(X),\Q)$
is an isomorphism and thus $g:=f\circ \iota$ has the desired property.
\end{proof}

\subsection{Moore spaces}
Given $G$ an abelian group, let $M(G,n)$, $n>1$, be the Moore space
(uniquely determined up to homotopy type) whose $n$-th homology group is $G$.
For example, $M(\Z,n)=S^n$ and $M(\Z_m,n)=e^{n+1}\cup_{f} S^{n}$, where
$f:S^n\to S^n$ is a map of degree $m$.
Note that $M(G,n)$ has the rational homotopy type of a wedge of spheres.

We will use the following properties of Moore spaces 
(cf. \cite[Proposition 1.7]{A}):

\begin{enumerate}
\item $M(A\oplus B,n)=M(A,n)\vee M(B,n)$;
\item A morphism $f:A \to B$  induces a continuous map $M(f):M(A,n)
\to M(B,n)$, so that $M(fg)=M(f)M(g)$ and $M(f)_* =f$. 
\end{enumerate}

In particular if $f:A\to B$ is an injection and there exists $g:B\to
A$ such that $gf=1_A$  then $M(B,n)$ dominates $M(A,n)$.
If $G$ is a finitely generated abelian group, then $M(G,n)$
is given by the wedge sum of copies of $S^n$ and copies of
$M(\Z_m,n)$ and $M(G,n)$ dominates any of these Moore subspaces.

\subsection{Homology decompositions}

Every simply connected $CW$ complex $Y$ has a {\it homology decomposition}
(cf. Theorem 4H.3 in \cite{hatcher} or \cite[Theorem 2.2]{A}). 
This means that there exists a homotopy equivalence $f:X\to Y$ such that $X$ can be constructed by the following iterated procedure.

Let $G_n:=H_n(Y)$. There exists an increasing sequence of complexes
$X_1\subset X_2\subset\cdots$ with $H_i(X_n)=G_n$ for $i\leq n$ and
$H_i(X_n)=0$ for $i>n$ where:

\begin{enumerate}

\item $X_1$ is a point and $X_2$ is the Moore space $M(G_2,2)$;

\item $X_{n+1}$ is the mapping cone of a cellular map
$h_n:M(G_{n+1},n)\to X_n$ such that the induced map
$(h_{n})_{*}:H_{n}(M(G_{n+1},n))\to H_n(X_n)$ is trivial;

\item $X=\cup_{n}X_{n}$.

\end{enumerate}

If $Y$ is a simply connected CW complex whose only non-zero homology
groups are $H_n(Y)$ and $H_{n+1}(Y)$, {\it and} $H_{n+1}(Y)$ is
a free abelian group, the homology
decomposition says in this case that $Y$ has the homotopy
type of $M(G_n,n)\vee M(G_{n+1},n+1)$ (cf. for example \cite[Lemma
2.6.5]{A}). It follows that:

\begin{Proposition} Suppose $Y$ is a simply connected CW complex whose only
non-zero homology groups are $H_n(Y)$ and $H_{n+1}(Y)$, where
$H_{n+1}(Y)$ is a free abelian group. Then $Y$ dominates
$M(H_n(Y),n)$ and $M(H_{n+1}(Y),n+1)$.
\label{util}
\end{Proposition}

\begin{Proposition}Suppose $X$ is a simply connected CW complex whose only
non-zero homology groups are $H_n(X)$ and $H_{n+1}(X)$.
\begin{enumerate}
\item if $b_{n}:=\mbox{\rm dim}\,H_{n}(X,\Q)\geq 2$ there
exists a map $g:S^n\vee S^n\to X$ such that
 $H_{*}(\Omega(g),\Q)$ is an injection;
\item if $H_{n+1}(X)$ is free abelian and different from $0$ or $\Z$, then
$X$ dominates $S^{n+1}\vee S^{n+1}$;
\item if $H_n(X)$ has an element of infinite order and $H_{n+1}(X)=\Z$, then
there exists a map $g:S^n\vee S^{n+1}\to X$ such that
 $H_{*}(\Omega(g),\Q)$ is an injection.
\end{enumerate}
\label{control}
\end{Proposition}

\begin{proof} Note that the Hurewicz map $\pi_{*}(X)\to H_{+}(X)$ is surjective.
Hence $X$ has the rational homotopy type of a wedge of spheres and by Lemma \ref{domrat},
if $b_n\geq 2$ there exists a map $g:S^n\vee S^n\to X$ such that
$H_{*}(\Omega(g),\Q)$ is an injection which proves the first item.
Similarly, if $H_n(X)$ has an element of infinite order and $H_{n+1}(X)=\Z$, then
Lemma \ref{domrat} gives a map $g:S^n\vee S^{n+1}\to X$ such that
$H_{*}(\Omega(g),\Q)$ is an injection, which proves the third item.

Finally, to prove the second item, note that since $H_{n+1}(X)$ is free
abelian and different from $0$ and $\Z$, there exists $\Z\oplus\Z\subset H_{n+1}(X)$ which is a direct summand. Hence using Proposition \ref{util}
we deduce that $X$ dominates $S^{n+1}\vee S^{n+1}$.

\end{proof}

\subsection{Ends of groups, $\ell_2$-Betti numbers and amenability}

The space of ends $E(X)$ of a locally compact separable metric space $X$
is given by the inverse limit
\[\lim_{K\subset X}\pi_{0}(X-K),\]
where the sets $K$ are compact. The space $E(X)$ is a totally disconnected
topological space and when $X$ is connected and locally connected, $E(X)$
is compact. Given a group $\pi$ acting freely  on a connected
simplicial complex $X$ with finite quotient, the homeomorphism type of
$E(X)$ only depends on $\pi$. The cardinality of $E(X)$ is usually denoted
by $e(\pi)$ and is called {\it number of ends of} $\pi$.
A finitely generated group $\pi$ has $0,1,2$ or infinitely many ends.
It has 0 ends if and only if it is finite. 

The higher order end groups of a group $\pi$
are defined as the cohomology groups of a $K(\pi,1)$ space with 
coefficients in the group ring $\Z[\pi]$.
We denote them by $H^{k}(\pi,\Z[\pi])$, $k\geq 0$.
If the group is infinite, then
$H^{1}(\pi,\Z[\pi])$ is a free abelian group of rank $e(\pi)-1$.

The group $\pi$ has two ends if and only if it is virtually $\Z$.
If $\pi$ has infinitely many ends, then it must contain a non-cyclic
free subgroup.

Finally, we note that if $\M$ is the universal covering of a closed
$n$-dimensional manifold $M$ with infinite fundamental group $\pi$, 
then $H^{1}(\pi,\Z[\pi])\cong H_{n-1}(\M)$.
Of course, one also has $H_{n}(\M)=0$.

Let $M$ be a closed manifold with an infinite amenable fundamental group
and let $\beta^{(2)}_{k}$ be the $k$-th $\ell_{2}$-Betti number of the universal
covering of $M$. It is interesting to note the following fact:

\begin{itemize}
\item If $\mbox{\rm dim}\,H_{k}(\M,\Q)$ is finite, then $\beta_{k}^{(2)}=0$ \cite{E1}.
\end{itemize}

This prompts the following question which is closely related 
to the Question in the introduction: 
if for some $k$, $\beta_{k}^{(2)}\neq 0$, is it true
that there exists a finite 1-connected rationally hyperbolic complex
$K$ and a map $f:K\to \M$ for which 
the rank of $H_{*}(\Omega(f),\Q)$ grows exponentially? 


We will make use of the following theorem.

\begin{Theorem} Let $M$ be a closed $n$-manifold with an infinite 
amenable fundamental group $\pi$ and let $\M$ be the universal covering
of $M$. Suppose $\M$ is $(k-1)$-connected, $k\geq 2$, and
$\mbox{\rm dim}\, H_{k}(\M,\Q)$ is finite. Then $H_{n-k}(\M,\Z)\cong H^{k}(\pi,\Z[\pi])$.
\label{eckext}
\end{Theorem}

\begin{proof} The proof is exactly the same as the proof of Theorem 3.1 in \cite{E2}.
The hypothesis $\mbox{\rm dim}\, H_{k}(\M,\Q)<\infty$ and the amenability of $\pi$
ensures that the $k$-th $\ell_{2}$-Betti number vanishes and one argues with
the commutative diagram on page 507 to conclude that 
$H_{\mbox{\rm comp}}^{k}(\M,\Z)\cong H^{k}(\pi,\Z[\pi])$. 
By Poincar\'e duality $H_{n-k}(\M,\Z)\cong H_{\mbox{\rm comp}}^{k}(\M,\Z)$.

\end{proof}

\section{Proof of Theorems A and B}

We first show:

\begin{Theorem} Let $M$ be a closed 4-manifold with infinite fundamental
group $\pi$. If $M$ admits a Riemannian metric with zero topological entropy,
then $\chi(M)=\sigma(M)=0$ and $H_{2}(\M)\cong H^{2}(\pi,\Z[\pi])$. 
Moreover, $\M$ has the rational homotopy type of a point, $S^2$ or $S^3$.
\label{main}
\end{Theorem}

\begin{proof} If $M$ admits a metric with zero entropy, $\pi$
has subexponential growth and hence it is amenable and can only have
1 or 2 ends.
Following B. Eckmann in \cite{E1}, we note that if $\pi$ is amenable
we can construct a F$\phi$lner sequence, that is, an increasing sequence
$Y_j$, $j=1,2,3,\cdots$, of finite subcomplexes of $\M$ with the 
following properties:
\begin{enumerate}
\item $Y_j$ consists of $N_j$ translates of a closed cellular fundamental
domain $D$ for the action of $\pi$;
\item $\cup_{j}Y_j=\M$;
\item let $\dot{N}_{j}$ be the number of translates of $D$ which meet the topological boundary of $Y_j$; then 
$$\lim_{j\to\infty}\frac{\dot{N}_{j}}{N_j}=0.$$
\end{enumerate}
Eckmann shows in \cite[p. 389]{E1} that
\[\chi(M)=\lim _{j\to\infty}\frac{b_2(Y_j)}{N_j},\]
where $b_{2}(Y_j)=\dim\,H_{2}(Y_j,\Q)$.
Moreover, by Proposition 2.1 in \cite{E1} we know that
if $b_{2}(\M)$ is finite, we must have
$\lim _{j\to\infty}\frac{b_2(Y_j)}{N_j}=0$.

Note that $H_1(\M)=H_4(\M)=0$ and hence the Hurewicz map $\pi_* (\M ) \to  H_{+} (\M )$ is onto 
and $\M$ has the rational homotopy type of a wedge of spheres

$$X:=\left(\vee_{\alpha}S^{2}_{\alpha}\right)\vee
\left(\vee_{\beta}S^{3}_{\beta}\right).$$
Since $X$ dominates any finite subcollection of them and 
the rational loop space homology of the wedge of at least two spheres 
grows exponentially, the Technical Lemma implies that $\M$ must
have the rational homotopy type of either a point, $S^2$ or $S^3$.
Thus $b_2(\M)$ is finite and $\chi(M)=0$ as desired.

To prove that $\sigma(M)=0$, we use the following observation
of Gromov in \cite[p. 85]{G2}: if $\sigma(M)\neq 0$, then $b_2(\M)$
must be infinite (this is a consequence of the amenability of $\pi$ and
the index theorem for infinite coverings).

Once we know $\chi(M)=0$, the isomorphism between $H_{2}(\M)$
and $H^{2}(\pi,\Z[\pi])$ is precisely Theorem 3.1 in \cite{E2} or 
Theorem \ref{eckext}.

\end{proof}

\begin{Remark}{\rm Let $M$ be a closed manifold of dimension $2k$ and suppose
that $\pi_i(M)=0$ for $1<i\leq k-1$ (the condition is vacuous
for $k=2$).
Suppose further that $\pi:=\pi_1(M)$ satisfies the following property:
it is infinite and the end groups $H^{i}(\pi,\Z[\pi])$ are zero
for $0<i<k$. It is quite easy to check (see Proposition 2.1 in \cite{E2})
that $H_{i}(\M)\cong H^{2k-i}(\pi,\Z[\pi])$ for $k<i\leq 2k$.
Thus $H_{k}(\M)$ is the only non-zero homology group.
If $M$ admits a metric with zero entropy, then $\chi(M)=0$.
The proof is the same as in the 4-dimensional case.
}
 \end{Remark}

\begin{Corollary} Let $M$ be a closed 4-manifold whose fundamental group
has two ends. If $M$ admits a metric with zero topological entropy,
it is finitely covered by $S^3\times S^1$.
\label{cor1}
\end{Corollary}

\begin{proof} Theorem 11.1 in \cite{H} says that a closed 4-manifold whose 
fundamental group has two ends and $\chi(M)=0$ is finitely covered
by $S^3\times S^1$.
\end{proof}

Finding the homeomorphism types of such manifolds is a fairly complicated
problem, we refer the interested reader to Chapter 11 in \cite{H}.

The last corollary and Theorem \ref{main} tell us that if we wish
to move further into the classification of closed 4-manifolds which admit
a metric of zero entropy we need to know more about $H^{2}(\pi,\Z[\pi])$
for $\pi$ with subexponential growth and one end. As far as we know, there is
no general result in this direction. 
However, note that if $\pi$ is the fundamental group of a closed manifold
whose universal covering is $\re^n$, then $H^{2}(\pi,\Z[\pi])$ is zero
if $n\neq 2$ and $\Z$ if $n=2$.

We can state:

\begin{Corollary} Let $M$ be a closed 4-manifold with $H^{2}(\pi,\Z[\pi])=\Z$.
If $M$ admits a metric with zero topological entropy, $M$ has a 
covering space of degree dividing 4 which is s-cobordant to $S^2\times \T^2$.
\label{cor2}
\end{Corollary}

\begin{proof} By Theorem \ref{main}, $\pi_{2}(M)\cong H_{2}(\M)\cong H^{2}(\pi,\Z[\pi])=\Z$. The corollary now follows from Theorem 10.1 in \cite{H}.
\end{proof}

There is no example known of a finitely {\it presented} group which is
of subexponential growth, but not of polynomial growth. Recall that 
the existence of a metric with zero entropy implies subexponential
growth of $\pi$.

\begin{Theorem}
Let $M$ be a closed 4-manifold whose fundamental group $\pi$
is infinite and has polynomial growth. 
If $M$ admits a metric with
zero topological entropy, then $M$ is finitely covered by one of the following:
\begin{enumerate}
\item $S^3\times S^1$;
\item a manifold s-cobordant to $S^2\times \T^2$;
\item a manifold homeomorphic
to a nilmanifold.
\end{enumerate}
\label{pol}

\end{Theorem}

\begin{proof} By a celebrated theorem of Gromov, $\pi$ is
virtually nilpotent. Thus by passing to a finite covering we can assume
that $\pi$ is nilpotent. It follows that $\pi$ coincides with its
Hirsch-Plotkin radical $\sqrt{\pi}$, which is the maximal nilpotent
normal subgroup.

Let $h(\pi)$ denote the Hirsch length of $\pi$. If $h(\pi)\leq 2$, then
up to finite index, $\pi$ must be $\Z$ or $\Z^2$.
If $\pi$ is $\Z$, then by Corollary \ref{cor1}, $M$ falls under item 1.
Similarly, if $\pi$ is $\Z^2$, by Corollary \ref{cor2}, $M$ falls
under item 2.

Finally if $h(\pi)=h(\sqrt{\pi})\geq 3$, Corollary 8.1.1 in \cite{H} implies
that $M$ is finitely covered by a manifold homeomorphic to a nilmanifold, since
by Theorem \ref{main}, $\chi(M)=0$.

\end{proof}

Theorems \ref{main} and \ref{pol} give Theorem A.

\subsection{Compact complex surfaces} Using the results on 4-manifolds
from the previous subsection, we can now complete the classification
of compact complex surfaces which admit a metric with
zero entropy. We begun this classification in \cite{PP}, but our results
excluded two cases:
\begin{enumerate}
\item surfaces of general type;
\item surfaces of class VII with positive second Betti number.
\end{enumerate}

It is well known that surfaces of general type have $\chi>0$.
A surface of class VII has first Betti number equal to one, thus a surface
of class VII with positive second Betti number must also have $\chi>0$.
Hence in both cases, Theorem \ref{main} implies that neither of these classes
admits a metric of zero entropy unless the fundamental group is finite.
Obviously, surfaces of class VII have an infinite fundamental group.
On the other hand, we know that a closed simply connected 4-manifold
that admits a metric with zero entropy must be homeomorphic to
to $S^{4}$, $\cp^{2}$, $S^{2}\times S^{2}$,
$\cp^{2}\#\overline{\cp}^{2}$ or $\cp^{2}\#\cp^{2}$ (cf. \cite{PP0}). 
Thus if there exists a surface of general type with a metric of
zero entropy and finite fundamental group, its universal covering
(which is also a surface of general type) would have to be homeomorphic
to $S^{2}\times S^{2}$ or $\cp^{2}\#\overline{\cp}^{2}$. It is unknown
whether there are such exotic examples, although
it is known that there is no surface of general type {\it diffeomorphic}
to $S^{2}\times S^{2}$ or $\cp^{2}\#\overline{\cp}^{2}$. 
Below we call such a potential example, {\it an exotic surface of general type}.

We now combine this discussion with the results in \cite{PP} to obtain:

\medskip

\noindent {\bf Theorem B.} {\it Let $S$ be a compact complex surface which is
not an exotic surface of general type.
Then $S$ admits a metric with zero topological entropy if and only if
$S$ is diffeomorphic to one of the following:
$\C P^2$, a ruled surface of genus $0$ or $1$, a complex torus, a hyperelliptic surface,
a Hopf surface, a Kodaira surface, or a Kodaira surface modulo a finite group.}

\section{Proof of Theorem C}

\noindent {\bf Theorem C.} {\it Let $M$ be a closed 5-manifold with infinite fundamental group. 
If $M$ admits a metric
with zero entropy, then $\M$ has the rational homotopy type of a finite 1-connected
elliptic complex. Moreover, $H_{3}(\M)\cong H^2(\pi,\Z[\pi])$.}

\medskip

\begin{proof}
We know that since $\pi$ is infinite amenable it can only have one or two ends.
Moreover, we know that $H_5(\M)=0$ and $H_{4}(\M)\cong H^{1}(\pi,\Z[\pi])$.

Suppose first that $\pi$ has one end. Since $H_4(\M)=0$, the Hurewicz map
$\pi_*(\M)\to H_{+}(\M)$ is surjective and $\M$ has the rational homotopy
type of a wedge of spheres 
$$X:=\left(\vee_{\alpha}S^{2}_{\alpha}\right)\vee \left(\vee_{\beta}S^{3}_{\beta}\right).$$
Since $X$ dominates and finite wedge $K$ of spheres from the collection
of $S^{2}_{\alpha}$ and $S^{3}_{\beta}$, we conclude that if 
$\mbox{\rm dim}\, H^{+}_{k}(\M,\Q)\geq 2$, there exist a finite simply connected
rationally hyperbolic complex $K$ and a map $f:K\to\M$ such that
$H_{*}(\Omega(f),\Q)$ is injective. Hence the Technical Lemma implies
that $\mbox{\rm dim}\, H^{+}_{k}(\M,\Q)\leq 1$ and $\M$ has the rational
homotopy type of a point, $S^2$ or $S^3$.

Suppose now that $\pi$ has two ends and so $H_4(\M)=\Z$.
Let us consider the homology decomposition of $\M$ and note that $\M$ 
has the homotopy type of the mapping cone of a cellular map $h_3:M(\Z,3)\to X_3$ such that
$(h_{3})_{*}:H_{3}(M(\Z,3))\to H_3(X_3)$ is trivial.
Since the Hurewicz map $\pi_*(X_3)\to H_{+}(X_3)$ is onto, $X_3$
has the rational homotopy type of a wedge of spheres
$$Y:=\left(\vee_{\alpha}S^{2}_{\alpha}\right)\vee \left(\vee_{\beta}S^{3}_{\beta}\right)$$
and let $f:X_3\to Y$ be a rational homotopy equivalence.
The image of the map $h_3\circ f$ can only intersect a finite number of
spheres from the collection. Therefore if the collection is infinite it
would follow that $\M$ has the rational homotopy type of a space
that dominates the wedge product of two
spheres. Since the rational loop space homology of the wedge of
two spheres grows exponentially the Technical Lemma would imply
positive entropy. Therefore there are only a finite number of
2-spheres
and 3-spheres in the collection and $\M$ has the rational
homotopy type of a finite CW complex (which of course must be
elliptic).

The conclusion $H_{3}(\M)\cong H^{2}(\pi,\Z[\pi])$ follows now 
directly from Theorem \ref{eckext}.

\end{proof}

\begin{Remark}{\rm The main point that makes the proof of Theorems \ref{main}
and C work is the fact that the rational Hurewicz map
$\pi_{*}(\M)\otimes\Q\to H_{+}(\M,\Q)$ is surjective. 
By Theorem 4.5 in \cite{FHT} this condition is equivalent to saying that
$\M$ has the rational homotopy type of a wedge of spheres or that
$\M$ has the rational homotopy type of a suspension.}
\end{Remark}

\subsection{Torsion of 5-manifolds}

\begin{Theorem} Let $M$ be a closed 5-manifold whose fundamental group
has one end and $H^{2}(\pi,\Z[\pi])$ is either $0$ or $\Z$.
If $M$ admits a metric with zero entropy, then:
\begin{enumerate}
\item $\mbox{\rm dim}\,H_2(\M,\Q)\leq 1$ and $H_2(\M)$ has no finite subgroup
as a direct summand; 
\item $H_{3}(\M)$ is either $0$ or $\Z$;
\item if $H_3(\M)=\Z$, then $H_2(\M)$ is a torsion group with no finite
subgroup as a direct summand.
\end{enumerate}
\end{Theorem}

An example of a group as in item 3 is the quasicyclic group of type
$p^{\infty}$ ($p$ prime) given by the $p^{k}$-th roots of unity, $k$ running
over all natural integers. We do not know if such a group can be realized
as $\pi_{2}(M)$ of a 5-manifold.

\begin{proof} Since the fundamental group of $M$ has one end, $H_{4}(\M)=0$.
By Theorem C, $H_3(\M)$ is $0$ or $\Z$ and we can apply
Proposition \ref{control} to $\M$. This proposition combined
with the Technical Lemma proves all the claims in the theorem
except the one regarding the abscence of finite groups of $H_2(\M)$
as a direct sumand. Suppose there is such a group. Then
there exists a prime $p$ such that the group $\Z_{p^k}$ appears
as a direct summand for some $k\geq 1$.
Then $\M$ dominates the Moore space $M(\Z_{p^{k}},2)$, but the latter
has the property that its loop space homology with coefficients in
$\Z_p$ grows exponentially. Again, this cannot happen by
the Technical Lemma.
\end{proof}


\begin{thebibliography}{aa}

\bibitem{AP} J.W. Anderson, G.P. Paternain, 
{\it The minimal entropy problem for 3-manifolds with zero simplicial volume}, 
Asterisque, {\bf 286} (2003) 63--79.


\bibitem{A} M. Aubry, {\it Homotopy theory and models,} based on lectures held at a DMV seminar in Blaubeuren by H. J. Baues, S. Halperin and
  J.-M. Lemaire. DMV Seminar, {\bf 24}, Birkhäuser Verlag, Basel, 1995.




\bibitem{D} E.I. Dinaburg, {\it On the relations among various entropy
characteristics of dynamical systems}, Math. USSR Izv. {\bf 5} (1971)
337--378.



\bibitem{E1} B. Eckmann, {\it Amenable groups and Euler characteristic,}
Comment. Math. Helvetici {\bf 67} (1992) 383--393.


\bibitem{E2} B. Eckmann, {\it Manifolds of even dimension
with amenable fundamental group}, Comment. Math. Helvetici
{\bf 69} (1994) 501--511.


\bibitem{FHT}Y.  F\'elix, S. Halperin, J.C. Thomas, {\it Rational homotopy theory,} Graduate Texts in Mathematics, {\bf 205}
  Springer-Verlag, New York, 2001. 


\bibitem{G1} M. Gromov, {\it Homotopical effects of dilatations,}  J. Diff. Geom. {\bf 13} (1978) 303--310.

\bibitem{G} M. Gromov, {\it Groups of polynomial growth and expanding maps,} 
 Inst. Hautes \'Etudes Sci. Publ. Math. {\bf 53} (1981) 53--73.


\bibitem{G2} M. Gromov, {\it Volume and bounded cohomology},
Publ. Math. IHES {\bf 56} (1982) 1--99.



\bibitem{bates} M. Gromov, {\it Metric structures for Riemannian and non-Riemannian spaces,} Progress in Mathematics,
{\bf 152}, S.M. Bates (translator) 1999.








\bibitem{hatcher} A. Hatcher, {\it Algebraic Topology}, Cambridge University Press, 2002.


\bibitem{H} J.A. Hillman, {\it Four-manifolds, geometries and knots,}
 Geometry and Topology Monographs, {\bf 5} Geometry and Topology Publications,
 Coventry, 2002.


\bibitem{HMR} P. Hilton, G. Mislin, J. Roitberg, 
{\it Localization of nilpotent groups and spaces,} 
North-Holland Mathematics Studies, {\bf 15} 1975.







\bibitem{Man} R. Ma\~{n}\'{e}, {\it On the topological entropy of geodesic flows},
J. Diff. Geom. {\bf 45} (1997) 74--93.


\bibitem{Ma} A. Manning, {\it Topological entropy for geodesic flows},
Ann. Math. {\bf 110} (1979), 567--573

\bibitem{P} G.P. Paternain, {\it Geodesic Flows},
Progress in Mathematics, Birkh\"auser {\bf 180} 1999.




\bibitem{PP0} G.P. Paternain, J. Petean,
{\it Minimal entropy and collapsing with curvature
bounded from below}, Invent. Math. {\bf 151} (2003) 415-450.



\bibitem{PP} G.P. Paternain, J. Petean, {\it Entropy and collapsing of
compact complex surfaces,} to appear in the Proceedings of the London Math. Soc.




\bibitem{Y} Y. Yomdin, {\it Volume growth and entropy}, Israel J. Math. {\bf 57} (1987) 287--300.









\end{thebibliography}
\end{document}